\documentclass[11pt]{amsart}
\usepackage{amssymb,amsmath,txfonts,mathrsfs,mathtools, titletoc, tikz, stmaryrd} 
\usepackage[alphabetic]{amsrefs}

\usepackage{graphicx}

\newcommand{\mres}{\mathbin{\vrule height 1.6ex depth 0pt width
0.13ex\vrule height 0.13ex depth 0pt width 1.3ex}}
\newtheorem{theorem}{Theorem}[section]
\newtheorem{prop}[theorem]{Proposition}
\newtheorem{lemma}[theorem]{Lemma}
\newtheorem{remark}[theorem]{Remark}
\newtheorem{claim}[theorem]{Claim}

\newtheorem{question}[theorem]{Question}
\newtheorem{definition}[theorem]{Definition}
\newtheorem{cor}[theorem]{Corollary}

\numberwithin{equation}{section}

\def\pf{{\it Proof:}~}

\begin{document}

\title[The $(p, m)$-width of Riemannian manifolds and its realization]{The $(p, m)$-width of Riemannian manifolds and its realization}

\author{Guoyi Xu}
\address{Yau Mathematical Sciences Center, Jin Chun Yuan West Building
\\Tsinghua University, Beijing\\P. R. China, 100084}
\email{gyxu@math.tsinghua.edu.cn}
\date{\today}

\begin{abstract}
While studying the existence of closed geodesics and minimal hypersurfaces in compact manifolds, the concept of width was introduced in different contexts. Generally, the width is realized by the energy of the closed geodesics or the volume of minimal hypersurfaces, which are found by the Minimax argument. Recently, Marques and Neves used the $p$-width to prove the existence of infinite many minimal hypersurfaces in compact manifolds with positive Ricci curvature. However, whether the $p$-width can be realized as the volume of minimal hypersurfaces is not known yet. We introduced the concept of the $(p, m)$-width which can be viewed as the stratification of the $p$-width, and proved that the $(p, m)$-width can be realized as the volume of minimal hypersurfaces with multiplicities.
\end{abstract}

%\keywords{$p$-sweepout, Min-Max theory, $p$-width} \subjclass[2010]{53A10, 49Q05} 
\thanks{The author was partially supported by NSFC 11401336}
\dedicatory{Dedicated to Professor Robert Gulliver on the occasion of his $70$th birthday}

\maketitle

\section{Introduction}\label{SECTION introduction}

In $1917$, Birkhoff \cite{Birk} proved the following classical theorem:
\begin{theorem}[Birkhoff]\label{thm Birk}
{There exists a nontrivial closed geodesic for any metric on $\mathbb{S}^2$.
}
\end{theorem}

Birkhoff's idea to prove Theorem \ref{thm Birk} is the minimax method. More concretely, let $W^{1, 2}$ be the space of $W^{1, 2}$ maps from $\mathbb{S}^1$ to $\mathbb{S}^2$, and $\Xi$ is the set of continuous maps $\hat{\sigma}: (\mathbb{S}^1, 0)\rightarrow (W^{1, 2}, \Lambda^0 \mathbb{S}^2)$, where $\mathbb{S}^1= [0, 1]/\partial [0, 1]$, $\Lambda^0 \mathbb{S}^2$ is the set of all the point curves in $\mathbb{S}^2$. For any map $\hat{\sigma}\in \Xi$, the \textbf{homotopy-energy-width} of $\hat{\sigma}$ is defined as the following: 
\begin{align}
\mathcal{HEW}(\hat{\sigma})= \inf_{\sigma\in \Xi_{\hat{\sigma}}} \max_{t\in \mathbb{S}^1} E\big(\sigma (t)\big) \nonumber 
\end{align}
where the energy $E\big(\sigma (t)\big)= \int_{\mathbb{S}^1} \big|\partial_x \sigma (t)\big|^2 dx$ is defined for any $\sigma(t)\in W^{1, 2}$, and the homotopy class $\Xi_{\hat{\sigma}}$ is the set of all maps $\sigma\in \Xi$ which are homotopic to $\hat{\sigma}$ through maps in $\Xi$. When $\Xi_{\hat{\sigma}}$ is a nontrivial homotopy class, the width $\mathcal{HEW}(\hat{\sigma})$ is positive and realized by the energy of a nontrivial closed geodesic. One key observation above is that the isomorphism between $\pi_1(W^{1, 2})$ and $\pi_{2}(\mathbb{S}^2)$, which guarantees us to find the nontrivial homotopy class $\Xi_{\hat{\sigma}}$ with the positive width.

In $1951$, by generalizing Birkhoff's idea, Lyusternik and Fet \cite{LF} proved the following general theorem about the existence of closed geodesics (also see \cite{Fet}).
\begin{theorem}[Lyusternik-Fet]\label{thm existence of closed geodesics on general manifolds}
{On every compact Riemannian manifold $M^{n+ 1}$ with $n\geq 1$, there exists a closed geodesic.
}
\end{theorem}

Their idea is to consider the set $\Xi^{k}$, which is the set of continuous maps $\hat{\sigma}: (\mathbb{S}^k, 0)\rightarrow (W^{1, 2}, \Lambda^0 M)$. For any map $\hat{\sigma}\in \Xi^k$, the homotopy-energy-width of $\hat{\sigma}$ is defined similarly as above: 
\begin{align}
\mathcal{HEW}(\hat{\sigma})= \inf_{\sigma\in \Xi^k_{\hat{\sigma}}} \max_{t\in \mathbb{S}^k} E\big(\sigma (t)\big) \nonumber 
\end{align}
where the homotopy class $\Xi^k_{\hat{\sigma}}$ is the set of all maps $\sigma\in \Xi^k$ which are homotopic to $\hat{\sigma}$ through maps in $\Xi^k$. If $\Xi^k_{\hat{\sigma}}$ is a nontrivial homotopy class, the width $\mathcal{HEW}(\hat{\sigma})$ is positive and realized by the energy of a nontrivial closed geodesic. Similar as $\mathbb{S}^2$ case, the isomorphism between $\pi_k(W^{1, 2})$ and $\pi_{k+ 1}(M)$, which guarantees us to find the nontrivial homotopy class $\Xi^k_{\hat{\sigma}}$ with the positive width for some $k$.

People tried to use the above method to show the existence of closed minimal hypersurfaces in compact manifolds $M^{n+ 1}$. In the rest of the paper, unless otherwise mentioned, $M$ is always $(n+ 1)$-dim close Riemannian manifold.

In $1962$, Almgren \cite{Alm} studied the relationship between the homotopy groups of the integral cycle groups on closed manifold $M$ and the homology groups of the manifold $M$. Specifically, he proved the following theorem:
\begin{theorem}[Almgren]\label{thm Alm 1962}
{For $(n+ 1)$-dim closed Riemannian manifold $(M^{n+ 1}, g)$, 
\begin{align}
\pi_k\big(\mathcal{Z}_n^0(M^{n+ 1}; \mathbb{Z}_2)\big)\cong H_{n+ k}(M^{n+ 1}; \mathbb{Z}_2)\ , \quad \quad \quad \quad \forall k\in \mathbb{N} \nonumber 
\end{align}
where $\mathcal{Z}_n^0(M; \mathbb{Z}_2)$ is the path connected component of integral cycle space $\mathcal{Z}_n(M; \mathbb{Z}_2)$ containing $\mathbf{0}$.
}
\end{theorem}

From Theorem \ref{thm Alm 1962}, we have 
\begin{align}
\pi_1\big(\mathcal{Z}_n^0(M^{n+ 1}; \mathbb{Z}_2); \mathbb{Z}_2\big)= H_{n+ 1}(M^{n+ 1}; \mathbb{Z}_2)= \mathbb{Z}_2\neq 0 \nonumber 
\end{align}
hence one can find a continuous map $\hat{\sigma}: \left(\mathbb{S}^1, 0\right)\rightarrow \left(\mathcal{Z}_n(M; \mathbb{Z}_2), \mathbf{0}\right)$, whose homotopy class $\Theta_{\hat{\sigma}}$ is nontrivial. 

Furthermore, Almgren \cite{Alm2} developed a discretization method to define the width. Very roughly, he discretized the continuous map $\hat{\sigma}$ to get a homotopy sequence of mapping $S= \{\varphi\}_{i= 1}^{\infty}$, where $\varphi_i$ maps the suitably chosen discrete points of $\mathbb{S}^1$ to $\mathcal{Z}_n(M; \mathbb{Z}_2)$. Similarly, the homotopy class $\Theta_{\hat{\sigma}}$ can also be discretized as $\Pi_{S}$. The \textbf{discrete homotopy width of $S$} can be defined as the following:
\begin{align}
\mathcal{DHW}(S)= \mathcal{DHW}\left(\{\varphi\}_{i= 1}^{\infty}\right)=
\inf_{\{\hat{\varphi}_i\}_{i= 1}^{\infty}\in \Pi_S}\left[ \varlimsup_{i\rightarrow \infty} \max_{x\in dmn(\hat{\varphi}_i)}\Big\{\mathbf{M}\big(\hat{\varphi}_i(x)\big)\Big\}\right] \nonumber
\end{align}
where $dmn(\hat{\varphi}_i)$ denotes the domain of $\hat{\varphi}_i$ and $\mathbf{M}(\cdot)$ is the mass of the current.

One advantage of discrete homotopy sequence over the continuous maps is that there is always so called critical sequence in the discrete homotopy class $\Pi_S$. Using this fact, Almgren \cite{Alm2} constructed a variational calculus in the large, from which he concluded that the width $\mathcal{DHW}(S)$ is positive and can be realized by the mass of a nonzero stationary integral $n$-varifold.

However, Almgren's theorem on the existence of stationary integral varifolds is inadequate to settle the question of existence of regular minimal hypersurfaces on manifolds. This is because integral varifolds which are only stationary have in general essential singularities, possibly of positive measure.

%Through much more complicated argument, Pitts got the result about minimal hypersurfaces similar as the closed geodesics case. More concretely, he proved that when $\Pi_{S}$ is a nontrivial homotopy class (note $\Pi_S$ is homotopically nontrivial if the corresponding $\Theta_{\hat{\sigma}}$ is homotopically nontrivial), the width $\mathcal{DHW}(\hat{\sigma})$ is positive and realized by the mass of a nontrivial close minimal hypersurface (with possible singularities). 

In $1981$, J. Pitts \cite{Pitts} constructed a similar variational calculus as \cite{Alm2}, but made important progress. He proved that the width $\mathcal{DHW}(S)$ can be realized by the mass of a nonzero stationary integral $n$-varifold with an additional variational property------$\mathbb{Z}_2$ almost minimizing, which implies the singular set is empty for hypersurfaces of $M^{n+ 1}$, with $2\leq n\leq 5$. Soon after Pitts' work, R. Schoen and L. Simon \cite{SS} generalized the results of Pitts' to $2\leq n\leq 6$ among a general context. We state their results in the form of the following theorem:
\begin{theorem}[Pitts, Schoen-Simon]\label{thm Pitts exist}
{If $S$ is a nontrivial $(I, \mathbf{M})$-homotopy sequence of mappings into $\left(\mathcal{Z}_n^0(M^{n+ 1}; \mathbf{M};\mathbb{Z}_2), \{\mathbf{0}\}\right)$, with $2\leq n\leq 6$, then 
\begin{align}
\mathcal{DHW}(S)= ||V||(M) \nonumber 
\end{align}
where $V\in \mathcal{IV}_n(M)$ is an integral varifold supported by a smooth, closed, embedded minimal hypersurface with possible multiplicities.
}
\end{theorem}

%One consequence of his result is the following existence theorem:
%\begin{theorem}[\cite{Pitts}, \cite{SS}]\label{thm Pitts}
%{Let $M^{n+ 1}$ be an $(n + 1)$-dimensional smooth closed Riemannian manifold, $n\geq 2$. There is an embedded minimal hypersurface $\Sigma\subset M$,  which is a nontrivial hypersurface without boundary and with a singular set $Sing(\Sigma)$ of Hausdorff dimension at most $(n- 7)$.
%}
%\end{theorem}

The width of the homotopy class has also been studied in mean curvature flow \cite{CM-wm}, and Ricci flow \cite{CM-wr1} and \cite{CM-wr2}.

Come back to the closed geodesic case, for different homotopy class, we possibly get different homotopy-energy-width and even different closed geodesics. To guarantee getting more than one geometrically different closed geodesics, Lyusternik and {\v{S}}nirel{\cprime}man \cite{LS} studied the width of the homology class of the integral $1$-cycle space $\mathcal{Z}_1(\mathbb{S}^2; \mathbb{Z}_2)$ (also see \cite{Jost}). More concretely, they define the \textbf{homology-energy-width} as the following:
\begin{align}
\mathcal{W}_i(\mathbb{S}^{2})= \inf_{\sigma\in \mathfrak{H}_i}\max_{t\in dmn(\sigma)} \mathbf{M}\big(\sigma(t)\big) \ , \quad \quad \quad \quad i= 1, 2, 3\nonumber 
\end{align}
where $\mathfrak{H}_i$ is the set of maps $\sigma: X\rightarrow \mathcal{Z}_1(\mathbb{S}^2; \mathbb{Z}_2)$ such that 
$[\sigma(X)]= \tau_i\in H_i(\mathcal{Z}_1(\mathbb{S}^2; \mathbb{Z}_2); \mathbb{Z}_2)$, $X$ is a cubical subcomplex of $I^m= [0, 1]^m$, $\tau_i$ is the chosen non-zero element of $H_i(\mathcal{Z}_1(\mathbb{S}^2; \mathbb{Z}_2); \mathbb{Z}_2)$ and $[\sigma(X)]$ is the representative class of $\sigma(X)$ in $H_i(\mathcal{Z}_1(\mathbb{S}^2; \mathbb{Z}_2); \mathbb{Z}_2)$. Using the homology-energy-width, they proved the following theorem:
\begin{theorem}[Lyusternik-{\v{S}}nirel{\cprime}man]\label{thm L-S}
{For any metric on $\mathbb{S}^2$, there are at least three simple closed geodesics.
}
\end{theorem}

From the view point of the waist inequality, starting at $1983$, Gromov \cite{Gromov83}, \cite{Gromov88} and \cite{Gromov03} studied the width of a homology class of mappings, and get the upper and lower bound of the $p$-width for integral currents in the unit Euclidean ball (see \cite[Section $8$]{Gromov03}). In $2009$, L. Guth \cite{Guth} refined, generalized the idea of Gromov, and gave a complete proof of the bound of the $p$-width for a large family of homology class of mappings into the unit Euclidean ball. It is worth mentioning that Lyusternik-{\v{S}}nirel{\cprime}man's idea was used to get the lower bound of the $p$-width.

On the other hand, in $1982$ S.-T. Yau \cite{Yau} posed the conjecture that every compact Riemannian three-manifold admits infinite number of smooth, closed, immersed minimal surfaces. 

In $2013$, to study Yau's above conjecture, inspired by the work of Lyusternik and {\v{S}}nirel{\cprime}man, Gromov, Guth, the $p$-sweepout and $p$-width for compact Riemannian manifolds were defined by F. Marques and A. Neves \cite{MN2} as the following. By Theorem \ref{thm Alm 1962} and the results of algebraic topology, we can get the singular homology and cohomology of $\mathcal{Z}_n^0(M^{n+ 1}; \mathbb{Z}_2)$: 
\begin{align}
H^p \big(\mathcal{Z}_n^0(M^{n+ 1}; \mathbb{Z}_2); \mathbb{Z}_2\big)= \mathbb{Z}_2 \ , \quad \quad \quad \forall p\in \mathbb{N} \nonumber 
\end{align}
Assume the generator of $H^p \big(\mathcal{Z}_n^0(M^{n+ 1}; \mathbb{Z}_2); \mathbb{Z}_2\big)= \mathbb{Z}_2$ is $h_p$, and $X$ is a cubical subcomplex of $I^{k}= [0, 1]^k$ for some $k\in \mathbb{N}$.

\begin{definition}\label{def p-sweepout and p-width}
{A continuous map $\Phi: X\rightarrow \mathcal{Z}_n^0(M^{n+ 1}; \mathbb{Z}_2)$ is a \textbf{$p$-sweepout} if $\Phi^*(h_p)\neq 0\in H^p(X; \mathbb{Z}_2)$. The \textbf{$p$-width} of $M^{n+ 1}$ is 
\begin{align}
\omega_p(M)= \inf_{\Phi\in \mathscr{P}_p} \max_{x\in dmn(\Phi)} \mathbf{M}(\Phi(x)) \nonumber
\end{align}
where $\mathscr{P}_p$ is the set of all $p$-sweepouts $\Phi$ that have no concentration of mass.
}
\end{definition}

Furthermore, they developed Lyusternik-{\v{S}}nirel{\cprime}man theory for the $p$-sweepout in Riemannian manifolds. Combining with the related $p$-width bounds, the following theorem was proved, which confirmed Yau's conjecture partially.
\begin{theorem}[Marques-Neves]
{In any compact Riemannian manifolds $(M^{n+ 1}, g)$ with positive Ricci curvature and $2\leq n\leq 6$, there are infinitely many smooth, closed, embedded minimal hypersurfaces.
}
\end{theorem}

%In \cite{MN1}, Marques and Neves showed that $\omega_5(\mathbb{S}^3)= 2\pi^2$?????, and used it to confirm the well-known Willmore conjecture.

We recall the min-max definition of the $p^{th}$-eigenvalue of $(M, g)$. Consider the Rayleigh quotient
\begin{align}
E: V\rightarrow [0, \infty)\ , \quad \quad \quad E(f)= \frac{\int_M |\nabla f|^2 dV_g}{\int_M f^2 dV_g} \nonumber 
\end{align}
where $V\vcentcolon = W^{1, 2}(M)- \{0\}$, then
\begin{align}
\lambda_p= \inf_{p-plane\ P\subset V} \max_{f\in P} E(f) \nonumber 
\end{align}
It is well-known that for each integer $p\geq 1$, there exists $p^{th}$-eigenfunction $f_{p}\in V$ such that $E(f_{p})= \lambda_p$. As proposed by Gromov \cite{Gromov88}, $\omega_p(M)$ can be thought as a nonlinear analogue of the Laplace spectrum of $M$. A natural question is:
\begin{question}\label{ques main que}
{If $2\leq n\leq 6$, does there exist a varifold $V$ of a smooth, closed, embedded minimal hypersurface in $M^{n+ 1}$, with possible multiplicities; such that $\omega_p(M)= \|V\|(M)$?
}
\end{question}

\begin{remark}\label{rem generic realization}
{Marques and Neves \cite{MN2} raised a similar question in terms of $p$-sweepout. Motivated by Uhlenbeck's classical result \cite{Uhlenbeck} that generic metrics have simple eigenvalues, it was conjectured in \cite{Neves} that for generic metrics at least, $\omega_p(M)$ is achieved by a multiplicity one minimal hypersurface with index $p$.
}
\end{remark}

In this paper, we introduce a refined concept the $(p, m)$-width, which can be viewed as the stratification of the $p$-width.
\begin{definition}\label{def RP-M-sw 0}
{If $X$ is a cubical subcomplex of the $m$-dimensional cube $I^m$, a continuous map $\Phi: X\rightarrow \mathcal{Z}_n^0(M^{n+ 1}; \mathbb{Z}_2)$ is a \textbf{$(p, m)$-sweepout} if $\Phi^*(h_p)\neq 0\in H_p(X; \mathbb{Z}_2)$. The set of all $(p, m)$-sweepouts that have no concentration of mass, is denoted by $\mathscr{P}_{p, m}$. The \textbf{$(p, m)$-width} of $M^{n+ 1}$ is defined by 
\begin{align}
\omega_{p, m}(M)= \inf_{\Phi\in \mathscr{P}_{p, m}}\max_{x\in dmn(\phi)}\mathbf{M}\big(\Phi(x)\big)\nonumber 
\end{align}
}
\end{definition}

\begin{remark}\label{rem the relation between p-width and (p, m)-width}
{It is obvious that 
\begin{align}
\lim_{m\rightarrow \infty} \omega_{p, m}(M)= \omega_p(M) \quad \quad and \quad \quad \omega_{p, m+ 1}(M)\leq \omega_{p, m}(M) \nonumber 
\end{align}
}
\end{remark}

The main result of this paper is the following realization theorem.
\begin{theorem}\label{thm main result}
{If $2\leq n\leq 6$, for each $m, p\in \mathbb{N}$ with $m\geq 2p+ 1$, there exists $V\in \mathcal{IV}_n(M^{n+ 1})$, where $V$ is the varifold of a smooth, closed, embedded minimal hypersurface, with possible multiplicities; such that $\omega_{p, m}(M)= ||V||(M)$.
}
\end{theorem}

We explain the origin of the idea behind our proof. In $1982$, F. Smith (under supervision of L. Simon) gave a powerful variant of Pitts' variational calculus approach in his unpublished thesis \cite{Smith} (a complete proof of the Simon-Smith Theorem was also provided by T. Colding and C. De Lellis \cite{CD}). Later, De Lellis and D. Tasnady \cite{DT} generalized the variant approach of Simon-Smith in $3$-dimensional manifolds to any $(n+ 1)$-dimensional manifolds ($n\geq 2$) case. 

More specifically, they define a \textbf{generalized smooth map} $\Gamma: [0, 1]\rightarrow \mathcal{S}(M^{n+ 1})$, where $\mathcal{S}(M^{n+ 1})$ is the set of all closed subsets of $M$ with finite $n$-dim Hausdorff measure (denoted by $\mathcal{H}^n$) and satisfies
\begin{itemize}
\item For each $t$ there is a finite set $P_t\subset M$ such that $\Gamma(t)$ is a smooth hypersurface in $M\backslash P_t$; $\mathcal{H}^n(\Gamma(t))$ depends smoothly on $t$ and $t\rightarrow \Gamma(t)$ is continuous in the Hausdorff sense; on any $U\subset \subset M\backslash P_{t_0}$, one have $\lim\limits_{t\rightarrow t_0}\Gamma(t)= \Gamma(t_0)$ smoothly in $U$. 
\end{itemize}
A generalized smooth map $\Gamma$ is a \textbf{sweepout} of $M$ if there exists a family of $\{\Omega_t\}_{t\in [0, 1]}$ of open sets, such that $(\Gamma(t)\backslash \partial \Omega_t)\subset P_t$ for any $t$, $\Omega_0= \emptyset$, $\Omega_1= M$ and 
\begin{align}
\lim_{t\rightarrow s} Vol(\Omega_t\backslash \Omega_s)+ Vol(\Omega_s\backslash \Omega_t)= 0 \nonumber 
\end{align}

The set of all sweepouts is denoted as $\mathsf{S}$. The homotopy class $\mathsf{S}_{\Gamma}$ is defined to be the set of all maps in $\mathsf{S}$ which are homotopic to $\Gamma$ through maps in $\mathsf{S}$.  A family $\Lambda$ of sweepouts is called \textbf{homotopically closed} if it contains the homotopy class of each of its elements. Similarly, one can define the width of homotopically closed family of sweepouts $\Lambda$ as the following:
\begin{align}
\mathcal{W} (\Lambda)= \inf_{\Gamma\in \Lambda}\Big[\max_{t\in [0, 1]}\mathcal{H}^n\left(\Gamma(t)\right)\Big] \nonumber 
\end{align}

We formulate their results in the form of the following theorem:
\begin{theorem}[Simon-Smith, Colding-DeLellis, DeLellis-Tasnady]\label{thm DT existence}
{For any homotopically closed family $\Lambda$ of sweepouts of $M^{n+ 1}$, then there exists an integral varifold $\Sigma\in \mathcal{IV}_n(M)$, which is supported by a smooth, closed, embedded minimal hypersurface (where multiplicity is allowed) with a singular set of Hausdorff dimension at most $(n- 7)$, such that $\mathcal{W} (\Lambda)= ||\Sigma||(M)$.
}
\end{theorem}

In the proof of the above theorem, the continuous pull-tight method to prove the existence of the stationary varifold and furthermore the almost minimizing varifold, was applied to the sequence of sweepouts of possibly different homotopy classes. Note every sweepout in their definition is specially $1$-sweepout in Definition \ref{def p-sweepout and p-width}.

By carefully checking the original work of Pitts \cite{Pitts}, we observed the discrete pull-tight method there can apply to, the sequence of the discretization of $p$-sweepout of different homotopy classes for any $p\in \mathbb{N}$. After this paper was circulated, Xin Zhou informed us the similar results for $1$-sweepout were obtained in his Ph.D thesis. 

In this paper, we replace the nontrivial homotopy sequence of mappings in Theorem \ref{thm Pitts exist} by the sequence of mappings in the same homology class (but possibly in different homotopy classes), and then apply the discrete pull-tight method in \cite{Pitts} to prove the existence of the stationary and almost minimizing varifold. Combining with the well-known regularity results about the stationary and almost minimizing varifold, we obtain the main result of this paper.

The organization of this paper is as the following. In Section $2$, we recalled the results needed in this paper from Almgren-Pitts Min-Max theory. In Section $3$, we introduced the concepts of the $(p, m)$-width and the $(p, m)$-$\mathbf{M}$-width, and proved the equivalence between them. Then we use this fact to show the equivalence between the $p$-width and the $p$-$\mathbf{M}$-width. One reason of this indirect proof is that the estimate constant $C_0$ in Proposition \ref{prop MN2 3.10} (iii) depends on the dimension of the ambient space $I^m$ of the domain $X$, i.e. $m$.

In Section $4$, we introduced the $(p, m)$-homology sequence to replace the role of the homotopy sequence, proved the existence of stationary varifolds by the classical pulling tight method. For the proof of the existence of almost minimizing varifolds with mass equal to the $(p, m)$-$\mathbf{M}$-width, we follow very closely the exposition of \cite[$4.8-4.10$]{Pitts}. One thing worth mentioning is that to use the combinatorics lemmas of \cite{Pitts}, we have to fix the value of $m$, that is the reason we can not prove the realization of the $p$-width by the current method. Finally, combining with the regularity theory about almost minimizing varifolds, also including the equivalence between the different concepts of width, the realization of the $(p, m)$-width was proved.

\section{Background on Almgren-Pitts Min-Max theory}\label{SECTION 2}

In this section, we review some basic facts in Geometric Measure Theory, and the results in the Almgren-Pitts Minimax theory which will be used in the later sections.

The space $\mathcal{I}_k(M; \mathbb{Z}_2)$ of $k$-dimensional mod $2$ integral currents in M (see \cite[$4.2.26$]{GMT} for more details); the group $\mathcal{Z}_k(M; \mathbb{Z}_2)$ of mod $2$ flat chains, contains $T\in \mathcal{I}_k(M; \mathbb{Z}_2)$ with $\partial T = 0$. For closed subset $B\subset M$, the space $\mathcal{Z}_k(M, B; \mathbb{Z}_2)= \big\{T\in \mathcal{I}_k(M; \mathbb{Z}_2)\big|\ spt(\partial T)\subset B\big\}$, where $spt(\partial T)$ is the support of $\partial T$ (see \cite[$4.4.6$]{GMT} for more details);

The closure $\mathcal{V}_k(M)$, in the weak topology, of the space of $k$-dimensional
rectifiable varifolds in $M$. The space of integral rectifiable $k$-dimensional varifolds
in $M$ is denoted by $\mathcal{IV}_k(M)$.

Given $T\in \mathcal{I}_k(M; \mathbb{Z}_2)$, we denote by $|T|$ and $\|T\|$ the integral varifold and the Radon measure in $M$ associated with $T$, respectively. Given $V\in \mathcal{V}_k(M)$, $\|V\|$ denotes the Radon measure in $M$ associated with $V$. 

The above spaces come with several relevant metrics. The flat metric
and the mass of $T\in \mathcal{I}_k(M; \mathbb{Z}_2)$, denoted by $\mathcal{F}(T)$ and $\mathbf{M}(T)$, are defined in \cite[page $423$]{GMT} and \cite[page $426$]{GMT}, respectively. The $\mathbf{F}$-metric on $\mathcal{V}_k(M)$ is
defined in \cite[page $66$]{Pitts} and induces the varifold weak topology
on $\mathcal{V}_k(M)$. Finally, the $\mathfrak{F}$-metric on $\mathcal{I}_k(M; \mathbb{Z}_2)$ is defined by
\begin{align}
\mathfrak{F}(S, T) = \mathcal{F}(S- T) + \mathbf{F}(|S|, |T|) \label{def of frak F}
\end{align}

\begin{lemma}\label{cor 4.1 of MN1}
{Let $\mathcal{S}$ be a compact subset of $\mathcal{Z}_n(M; \mathfrak{F}; \mathbb{Z}_2)$, for every $\epsilon> 0$, there is $\delta> 0$ so that for every $S\in \mathcal{S}$ and $T\in \mathcal{Z}_n(M; \mathbb{Z}_2)$, 
\begin{align}
\mathbf{M}(T)< \mathbf{M}(S)+ \delta \quad and \quad \mathcal{F}(T- S)\leq \delta \Longrightarrow \mathfrak{F}(S, T)\leq \epsilon \nonumber 
\end{align}
}
\end{lemma}

\pf
{See \cite[Lemma $4.1$]{MN1}.
}
\qed

We assume that $\mathcal{I}_k(M; \mathbb{Z}_2)$, $\mathcal{Z}_k(M; \mathbb{Z}_2)$, have the topology induced by the flat metric. When endowed with the topology of the
mass norm, these spaces will be denoted by $\mathcal{I}_k(M; \mathbf{M}; \mathbb{Z}_2)$, $\mathcal{Z}_k(M; \mathbf{M}; \mathbb{Z}_2)$, respectively. The space $\mathcal{V}_k(M)$ is considered with the weak topology of varifolds. 

For each $j\in \mathbb{N}$, $I(1, j)$ denotes the cell complex of the unit interval $I^1= [0, 1]$ whose $1$-cells and $0$-cells are, respectively,
\begin{align}
[0, 3^{-j}],\ [3^{-j}, 2\cdot 3^{-j}], \cdots , [1- 3^{-j}, 1] \quad and \quad [0],\ [3^{-j}], \cdots , [1- 3^{-j}], [1] \nonumber 
\end{align}

We denote by $I(m, j)$ the cell complex on $I^m$:
\begin{align}
I(m, j)= I(1, j)\otimes\cdots \otimes I(1, j) \quad \quad \quad (m\ times) \nonumber 
\end{align}

Let $X$ be a cubical subcomplex of the $m$-dimensional cube $I^m= [0, 1]^m$. The cube complex $X(j)$ is the union of all cells of $I(m, j)$ whose support is contained in some cell of $X$, and $X(j)_q$ denotes the set of all $q$-cells in $X(j)$.

Given $i, j\in \mathbb{N}$, we define 
\begin{align}
\mathbf{n}(i, j): X(i)_0\rightarrow X(j)_0 \nonumber  
\end{align}
so that $\mathbf{n}(i, j)(x)$ is the element in $X(j)_0$ that is closest to $x$.

Two vertices $x, y\in X(j)_0$ are adjacent if they belong to a common cell in $X(j)_1$. Given a map $\varphi: X(j)_0\rightarrow \mathcal{Z}_n(M; \mathbf{M}; \mathbb{Z}_2)$, we define the \textbf{fineness of $\varphi$} to be 
\begin{align}
\mathbf{f}(\varphi)= \sup\{\mathbf{M}\big(\varphi(x)- \varphi(y)\big): x,\ y \ adjacent \ vertices\ in\ X(j)_0\} \nonumber
\end{align}

\begin{definition}\label{def p-homot}
{Assume $\varphi_j: X(k_j)_0\rightarrow \mathcal{Z}_n(M; \mathbf{M}; \mathbb{Z}_2)$, where $j= 1, 2$. We say \textbf{$\varphi_1$ is $X$-homotopic to $\varphi_2$ in $\mathcal{Z}_n(M; \mathbf{M}; \mathbb{Z}_2)$ with fineness $\delta$} if and only if there exists a positive integer $k_3$ and a map
\begin{align}
\psi : I(1, k_3)_0\times X(k_3)_0 \rightarrow \mathcal{Z}_n(M; \mathbf{M}; \mathbb{Z}_2) \nonumber
\end{align}
such that $\mathbf{f}(\psi)< \delta$,  and $x\in X(k_3)_0$, 
\begin{align}
\psi\big([j- 1], x\big)= \varphi_j\big(\mathbf{n}(k_3, k_j)(x)\big)\ , \quad \quad \quad \quad j= 1, 2\nonumber
\end{align}
}
\end{definition} 

\begin{definition}\label{def p-homot seq}
{An \textbf{$(X, \mathbf{M})$-homotopy sequence of mappings into $\mathcal{Z}_n(M; \mathbf{M}; \mathbb{Z}_2)$} is a sequence of mappings $S= \{\varphi_i\}_{i= 1}^{\infty}$,
\begin{align}
\varphi_i : X(k_i)_0\rightarrow \mathcal{Z}_n(M; \mathbf{M}; \mathbb{Z}_2) \nonumber 
\end{align}
such that $\varphi_i$ is $X$-homotopic to $\varphi_{i+ 1}$ in $\mathcal{Z}_n(M; \mathbf{M}; \mathbb{Z}_2)$ with fineness $\delta_i\rightarrow 0$ and $\sup\{\mathbf{M}(\varphi_i(x)): x\in X(k_i)_0,\ i\in \mathbb{N}\}< \infty$.
}
\end{definition}

\begin{definition}\label{def p-homot seq class}
{Let $S_1= \{\varphi_i^1\}_{i= 1}^{\infty}$ and $S_2= \{\varphi_i^2\}_{i= 1}^{\infty} $ be $(X, \mathbf{M})$-homotopy sequences of mappings into $\mathcal{Z}_n(M; \mathbf{M}; \mathbb{Z}_2)$, \textbf{$S_1$ is homotopic to $S_2$} if there exists a sequence $\delta_i\rightarrow \infty$, such that
$\varphi_i^1$ is $X$-homotopic to $\varphi_i^2$ in $\mathcal{Z}_n(M; \mathbf{M}; \mathbb{Z}_2)$ with fineness $\delta_i$.
}
\end{definition}

We denote by $\big[X, \mathcal{Z}_n(M; \mathbf{M}; \mathbb{Z}_2)\big]$ the set of all $(X, \mathbf{M})$-homotopy sequences of mappings into $\mathcal{Z}_n(M; \mathbf{M}; \mathbb{Z}_2)$. And let $dmn(\varphi_i)$ denote the domain of $\varphi_i$, we define $\mathbf{L}: \big[X, \mathcal{Z}_n(M; \mathbf{M}; \mathbb{Z}_2)\big]\rightarrow [0, +\infty]$ by setting 
\begin{align}
\mathbf{L}(S)= \varlimsup_{i\rightarrow \infty} \max_{x\in dmn(\varphi_i)} \mathbf{M}\big(\varphi_i(x)\big) \ , \quad \quad \forall S= \big\{\varphi_i\big\}_{i= 1}^{\infty}\in \big[X, \mathcal{Z}_n(M; \mathbf{M}; \mathbb{Z}_2)\big] \nonumber 
\end{align} 

\begin{definition}\label{def no mass concentration}
{Given a continuous map $\Phi: X \rightarrow \mathcal{Z}_n(M; \mathbb{Z}_2)$, we define
\begin{align}
\mathbf{m}(\Phi, r)= \sup_{x\in X, p\in M} \|\Phi(x)\|\big(B_p(r)\big) \nonumber 
\end{align}
We say that \textbf{$\Phi$ has no concentration of mass} if $\varlimsup\limits_{r\rightarrow 0}\mathbf{m}(\Phi, r)= 0$. 
}
\end{definition}

Define $X_b$ to be the boundary of $X$.

\begin{prop}\label{prop MN1 13.1}
{Let $\Phi: X\rightarrow \mathcal{Z}_n(M; \mathbb{Z}_2)$ be a continuous map in the flat topology that has no concentration of mass. There exists a sequence of maps
\begin{align}
\phi_i&: X(k_i)_0\rightarrow \mathcal{Z}_n(M; \mathbf{M}; \mathbb{Z}_2) \nonumber 
\end{align}
with $k_i< k_{i+ 1}$, and a sequence of positive number $\{\delta_i\}_{i\in \mathbb{N}}$ converging to $0$ such that 
\begin{enumerate}
\item[(a)] \begin{align}
\sup_{x\in X(k_i)_0} \mathbf{M}\big(\phi_i(x)\big) \leq \sup_{x\in X} \mathbf{M}\big(\Phi(x)\big)+ \delta_i  \nonumber
\end{align}

\item[(b)] $S= \{\phi_i\}_{i\in \mathbb{N}}$ is an $(X, \mathbf{M})$-homotopy sequence of mappings into $\mathcal{Z}_n(M; \mathbf{M}; \mathbb{Z}_2)$ with $\mathbf{f}(\phi_i)< \delta_i$.

\item[(c)] \begin{align}
\sup_{x\in X(k_i)_0}\mathcal{F}\big(\phi_i(x)- \Phi(x)\big)\leq \delta_i \nonumber 
\end{align}

\item[(d)] Moreover, if $\Phi\big|_{X_b}$ is continuous in the $\mathfrak{F}$-metric, then
\begin{align}
\mathbf{M}\big(\phi_i(x)\big)\leq \mathbf{M}\big(\Phi(x)\big)+ \delta_i\ , \quad \quad \quad \quad \forall x\in X_b(k_i)_0 \nonumber 
\end{align}
and if $\Phi\big|_{X_b}$ is continuous in the mass topology, one can choose $\phi_i$ so that 
\begin{align}
\phi_i(x)= \Phi(x) \ , \quad \quad \quad \quad \forall x\in X_b \nonumber 
\end{align}
\end{enumerate}
}
\end{prop}

\begin{remark}\label{rem discreitzation associated to continuous map}
{For any $\Phi$ and $\{\phi_i\}_{i= 1}^{\infty}$ satisfying $(a), (b), (c)$ as in Proposition \ref{prop MN1 13.1}, it is denoted as $\{\phi_i\}_{i= 1}^{\infty}= \mathscr{D}(\Phi)$, and $\{\phi_i\}_{i= 1}^{\infty}$ is called the \textbf{sequence of discretizations} associated to $\Phi$. 
}
\end{remark}

\pf
{See \cite[Theorem $13.1$]{MN1}.
}
\qed

\begin{prop}\label{prop MN2 3.10}
{There exist positive constants $C_0= C_0(M, m)$ and $\delta_0= \delta_0(M)$ so that if $Y$ is a cubical subcomplex of $I(m, k)$ and 
\begin{align}
\phi: Y_0\rightarrow \mathcal{Z}_n(M; \mathbf{M}; \mathbb{Z}_2) \nonumber 
\end{align}
has $\mathbf{f}(\phi)< \delta_0$, then there exists a map
\begin{align}
\Phi: Y\rightarrow \mathcal{Z}_n(M; \mathbf{M}; \mathbb{Z}_2) \nonumber 
\end{align}
continuous in the mass norm and satisfying
\begin{enumerate}
\item[(i)] $\Phi(x)= \phi(x)$ for all $x\in Y_0$;
\item[(ii)] If $\alpha$ is some $j$-cell in $Y_j$, then $\Phi$ restricted to $\alpha$ depends only on the values of $\phi$ assumed on the vertices of $\alpha$;
\item[(iii)] $\max\Big\{ \mathbf{M}\big(\Phi(x)- \Phi(y)\big): x, y$ lie in a common cell of $Y\Big\}\leq C_0 \mathbf{f}(\phi)$.
\end{enumerate}
}
\end{prop}

\begin{remark}\label{rem Alm exten}
{We call the map $\Phi$ given by Proposition \ref{prop MN2 3.10} the \textbf{Almgren extension} of $\phi$, and it is denoted as $\Phi= \mathscr{A}(\phi)$.
}
\end{remark}

\pf
{See \cite[Theorem $3.10$]{MN2} and \cite[Theorem $14.2$]{MN1}.
}
\qed

\begin{cor}\label{cor MN2 3.12}
{Let $S=\{\varphi_i\}_{i\in \mathbb{N}}$ and $S'= \{\varphi'_i\}_{i\in \mathbb{N}}$ be $(X, \mathbf{M})$-homotopy sequence of mappings into $\mathcal{Z}_n(M; \mathbf{M}; \mathbb{Z}_2)$ such that $S$ is homotopic with $S'$.
\begin{enumerate}
\item[(i)] If $\Phi_i= \mathscr{A}(\varphi_i)$ and $\Phi_i'= \mathscr{A}(\varphi'_i)$, then $\Phi_i$ is homotopic to $\Phi'_i$ in the flat topology for sufficiently large $i$.
\item[(ii)] If $S= \mathscr{D}(\Phi)$, where $\Phi: X\rightarrow \mathcal{Z}_n(M;\mathbb{Z}_2)$ is a continuous map with no concentration of mass, then $\Phi_i$ is homotopic to $\Phi$ in the flat topology for every sufficiently large $i$. Moreover,
\begin{align}
\varlimsup_{i\rightarrow \infty} \max_{x\in X} \mathbf{M}\left(\Phi_i(x)\right)= \mathbf{L}(S)\leq \max_{x\in X} \mathbf{M}\left(\Phi(x)\right) \nonumber 
\end{align}
\end{enumerate}
}
\end{cor}

\pf
{See \cite[Corollary $3.12$]{MN2}.
}
\qed

\section{The equivalence between the $p$-width and the $p$-$\mathbf{M}$-width}\label{SECTION 3}

In this section, we define the $(p, m)$-width and the $(p, m)$-$\mathbf{M}$-width, and proved the equivalence between them by Almgren's interpolation results. The equivalence between the $p$-width and the $p$-$\mathbf{M}$-width follows as the corollary of the above refined equivalence.

Denote $\mathcal{Z}_n^0(M; \mathbf{M}; \mathbb{Z}_2)$ as the path connected component of $\mathcal{Z}_n(M; \mathbf{M}; \mathbb{Z}_2)$ containing $\mathbf{0}$, similarly denote $\mathcal{Z}_n^0(M; \mathbb{Z}_2)$ as the path connected component of $\mathcal{Z}_n(M; \mathbb{Z}_2)$ containing $\mathbf{0}$. It is easy to see that $\mathcal{Z}_n^0(M; \mathbf{M}; \mathbb{Z}_2)\subset \mathcal{Z}_n^0(M; \mathbb{Z}_2)$.

\begin{lemma}\label{lem cohomology of cycle spaces}
{For any $p\in \mathbb{N}$, 
\begin{align}
H^1\big(\mathcal{Z}_n^0(M; \mathbb{Z}_2); \mathbb{Z}_2\big)= \langle \bar{\lambda} \rangle\ ,  \quad \quad \quad 
H^p\big(\mathcal{Z}_n^0(M; \mathbb{Z}_2); \mathbb{Z}_2\big)= \langle \bar{\lambda}^p\rangle = \mathbb{Z}_2 \nonumber 
\end{align}  
}
\end{lemma}

\pf
{From \cite{Alm} and \cite[Theorem $4.6$]{Pitts}, we have 
\begin{align}
\pi_k\left(\mathcal{Z}_n^0\left(M; \mathbb{Z}_2\right)\right)= H_{n+ k}(M^{n+ 1}; \mathbb{Z}_2)\ , \quad \quad \quad \forall k\geq 1 \nonumber 
\end{align} 
then by $H_{n+ 1}(M^{n+ 1}; \mathbb{Z}_2)= \mathbb{Z}_2$ and $H_{n+ k}(M^{n+ 1}; \mathbb{Z}_2)= 0$ when $k> 1$, we obtain
\begin{equation}\nonumber
\pi_k\left(\mathcal{Z}_n^0\left(M; \mathbb{Z}_2\right)\right)= \left\{
\begin{array}{rl}
&\mathbb{Z}_2 \quad  \quad \quad \quad \quad \quad \quad \quad k= 1\\
&0 \quad  \quad \quad \quad \quad \quad \quad \quad k= 0, 2, 3, \cdots\\
\end{array} \right.
\end{equation}

By \cite[Proposition $4.13$]{Hatcher}, there exists a CW-complex $Z$ and a weak homotopy equivalence $f: Z\rightarrow \mathcal{Z}_n^0\left(M; \mathbb{Z}_2\right)$, then 
\begin{equation}\nonumber
\pi_k\left(Z\right)= \left\{
\begin{array}{rl}
&\mathbb{Z}_2 \quad  \quad \quad \quad \quad \quad \quad \quad k= 1\\
&0 \quad  \quad \quad \quad \quad \quad \quad \quad k= 0, 2, 3, \cdots\\
\end{array} \right.
\end{equation}

Note the universal cover of $Z$, denoted as $\tilde{Z}$, is a CW-complex with $\pi_i(\tilde{Z})= 0$ for any nonnegative integer $i$. From \cite[Theorem $4.5$]{Hatcher}, we know that $\tilde{Z}$ is contractible, hence $Z$ is a Eilenberg-Maclane space $K(\mathbb{Z}_2, 1)$ defined as in \cite[$1.B$]{Hatcher}. It is well-known that $\mathbb{RP}^{\infty}$ is a $K(\mathbb{Z}_2, 1)$. By the uniqueness of the homotopic type of CW complex $K(\mathbb{Z}_2, 1)$, $Z$ is homotopic to $\mathbb{RP}^{\infty}$. So we have 
\begin{align}
H_p\big(Z; \mathbb{Z}_2\big)= H^p\big(Z; \mathbb{Z}_2\big)= H^p\big(\mathbb{RP}^{\infty}; \mathbb{Z}_2\big)= \mathbb{Z}_2\ , \quad \quad \quad \forall p\in \mathbb{N} \label{tp 1}
\end{align}

Because weak homotopy equivalence induces isomorphism between homology \cite[Proposition $4.21$]{Hatcher}, we get 
\begin{align}
H_p\big(\mathcal{Z}_n^0(M; \mathbb{Z}_2); \mathbb{Z}_2\big)= H_p\big(Z; \mathbb{Z}_2\big) \label{tp 2} 
\end{align}

From (\ref{tp 1}) and (\ref{tp 2}), we have 
\begin{align}
H_p\big(\mathcal{Z}_n^0(M; \mathbb{Z}_2); \mathbb{Z}_2\big)= \mathbb{Z}_2\ , \quad \quad \quad \forall p\in \mathbb{N} \nonumber 
\end{align}

From the Universal Coefficient Theorem in algebraic topology, we can assume 
\begin{align}
H^p\big(\mathcal{Z}_n^0(M; \mathbb{Z}_2); \mathbb{Z}_2\big)= \langle h_p\rangle= \mathbb{Z}_2 \nonumber 
\end{align} 

Hence we have 
\begin{align}
H^1\big(\mathcal{Z}_n^0(M; \mathbb{Z}_2); \mathbb{Z}_2\big)= \langle \bar{\lambda} \rangle \nonumber 
\end{align}
Denote by $\bar{\lambda}^p$ the cup product of $\bar{\lambda}$ with itself $p$ times, we get all of our conclusions.
}
\qed

We recall that $\mathscr{P}_{p, m}$ is the set of all continuous maps $\Phi: X\rightarrow \mathcal{Z}_n^0(M; \mathbb{Z}_2)$ having no concentration of mass, with $\Phi^*(\bar{\lambda}^p)\neq 0\in H^p(X; \mathbb{Z}_2)$, where $X$ is a cubical subcomplex of $I^m$. Note$\mathscr{P}_{p}= \bigcup_{m= 1}^{\infty} \mathscr{P}_{p, m}$.

We define $\mathscr{P}_{p, m}^{\mathbf{M}}$ as the set of all continuous maps $\Phi: X\rightarrow \mathcal{Z}_n^0(M; \mathbf{M}; \mathbb{Z}_2)$, such that $(\mathbf{Ic}\circ \Phi)^*(\bar{\lambda}^p)\neq 0$, where $X$ is a cubical subcomplex of $I^m$ and $\mathbf{Ic}: \mathcal{Z}_n^0(M; \mathbf{M}; \mathbb{Z}_2)\rightarrow \mathcal{Z}_n^0(M; \mathbb{Z}_2)$ is the continuous inclusion map.

\begin{definition}\label{def (p, m)-width}
{We define the \textbf{$(p, m)$-width} $\omega_{p, m}(M)$ and the \textbf{$(p, m)$-$\mathbf{M}$-width} $\omega^{\mathbf{M}}_{p, m}(M)$ as the following:
\begin{align}
\omega_{p, m}(M)&= \inf_{\Phi\in \mathscr{P}_{p, m}} \max_{x\in dmn(\Phi)} \mathbf{M}\left(\Phi(x)\right) \nonumber  \\
\omega^{\mathbf{M}}_{p, m}(M)&= \inf_{\Phi\in \mathscr{P}_{p, m}^{\mathbf{M}}} \max_{x\in dmn(\Phi)} \mathbf{M}\left(\Phi(x)\right) \nonumber  
\end{align} 
and if $\mathscr{P}_{p, m}= \emptyset$ ($\mathscr{P}_{p, m}^{\mathbf{M}}= \emptyset$), define $\omega_{p, m}= \infty$ ($\omega_{p, m}^{\mathbf{M}}= \infty$). 
}
\end{definition}

Then one can relate the $p$-width with the $(p, m)$-width as the following: 
\begin{align}
\omega_p(M)=\inf_{m\in \mathbb{N}}\omega_{p, m}(M) \label{def p-width 1}
\end{align}
And we define the \textbf{$p$-$\mathbf{M}$-width} of $M$ as:
\begin{align}
\omega_p^{\mathbf{M}}(M)= \inf_{m\in \mathbb{N}} \omega_{p, m}^{\mathbf{M}}(M) \label{def p-M-width} 
\end{align}
which will be shown to be equivalent to the $p$-width later.

We define $\mathscr{D}_{p, m}$ as the set of $S=\{\phi_i\}_{i\in \mathbb{N}}$, where $S$ is an $(X, \mathbf{M})$-homotopy sequence of mappings into $\mathcal{Z}_n^0(M; \mathbf{M}; \mathbb{Z}_2)$, with $X$ as a cubical subcomplex of $I^m$,  
\begin{align}
&\phi_i: X(k_i)_0\rightarrow \mathcal{Z}_n^0(M; \mathbf{M}; \mathbb{Z}_2) \nonumber \\
&k_i< k_{i+ 1}\ , \quad \mathbf{f}(\phi_i)< \delta_i\ , \quad \lim_{i\rightarrow \infty} \delta_i= 0 \nonumber 
\end{align} 
and $\Phi_i= \mathscr{A}(\phi_i)\in \mathscr{P}_{p, m}$ for sufficiently large $i$.

\begin{lemma}\label{lem MN2's lem 4.7}
{For any $p, m\in \mathbb{N}$, we have 
\begin{align}
\omega_{p, m}(M)= \inf\limits_{S\in \mathscr{D}_{p, m}}\mathbf{L}(S) \nonumber 
\end{align}
}
\end{lemma}

\pf
{Choose $S_{\epsilon}= \{\phi_i\}_{i\in \mathbb{N}}\in \mathscr{D}_{p, m}$ with 
\begin{align}
\mathbf{L}(S_{\epsilon})\leq \inf_{S\in \mathscr{D}_{p, m}} \mathbf{L}(S)+ \epsilon \label{mn2 4.7.1}
\end{align}
Let $\Phi_i= \mathscr{A}(\phi_i)\in \mathscr{P}_{p, m}$, by Proposition \ref{prop MN2 3.10} (i) and (iii), $\mathbf{f}(\phi_i)= \delta_i$ and $\lim\limits_{i\rightarrow \infty} \delta_i= 0$, 
\begin{align}
\omega_{p, m}(M)\leq \varlimsup_{i\rightarrow \infty} \max_{x\in dmn(\Phi_i)} \mathbf{M}\big(\Phi_i(x)\big) = \mathbf{L}(S_{\epsilon}) \label{mn2 4.7.2}
\end{align}
From (\ref{mn2 4.7.1}) and (\ref{mn2 4.7.2}), let $\epsilon\rightarrow 0$, we obtain
\begin{align}
\omega_{p, m}(M)\leq \inf_{S\in \mathscr{D}_{p, m}} \mathbf{L}(S) \label{mn2 4.7.3}
\end{align}

For any $\epsilon> 0$, choose $\Phi\in \mathscr{P}_{p, m}$, with
\begin{align}
\left[\max_{x\in dmn(\Phi)} \mathbf{M}\big(\Phi(x)\big) \right]\leq \omega_{p, m}(M)+ \epsilon \label{mn2 4.7.4} 
\end{align}

Consider $\tilde{S}= \{\phi_i\}_{i\in \mathbb{N}}= \mathscr{D}(\Phi)$, then from Proposition \ref{prop MN2 3.10} and \cite[Lemma $3.8$]{MN2}, $\tilde{S}\in \mathscr{D}_{p, m}$. By Corollary \ref{cor MN2 3.12} (ii), 
\begin{align}
\mathbf{L}(\tilde{S})\leq \max_{x\in dmn(\Phi)} \mathbf{M}\big(\Phi(x)\big) \label{mn2 4.7.5}
\end{align}

From (\ref{mn2 4.7.4}) and (\ref{mn2 4.7.5}), let $\epsilon\rightarrow 0$, we get 
\begin{align}
\inf_{S\in \mathscr{D}_{p, m}} \mathbf{L}(S)\leq \omega_{p, m}(M) \label{mn2 4.7.6}
\end{align}

By (\ref{mn2 4.7.3}) and (\ref{mn2 4.7.6}), the conclusion follows. 
}
\qed

\begin{prop}\label{prop width equation}
{$\omega_{p, m}^{\mathbf{M}}(M)= \omega_{p, m}(M)$ and $\omega_p^{\mathbf{M}}(M)= \omega_p(M)$. 
}
\end{prop}

\pf
{By \cite[Lemma $3.8$]{MN2} and the definition of $\mathscr{P}_{p, m}$, we get $\mathscr{P}_{p, m}^{\mathbf{M}}\subset \mathscr{P}_{p, m}$, 
\begin{align}
\omega_{p, m}^{\mathbf{M}}(M)\geq \omega_{p, m}(M) \label{equi 1.1}
\end{align}

From Lemma \ref{lem MN2's lem 4.7}, for any $\epsilon> 0$, there exists $S_{\epsilon}=\{\phi_i\}_{i\in \mathbb{N}}\in \mathscr{D}_{p, m}$, such that 
\begin{align}
\omega_{p, m}(M)+ \epsilon\geq \mathbf{L}(S_{\epsilon}) \label{equi 1.2}
\end{align}

For all sufficiently large $i$, let $\Phi_i= \mathscr{A}(\phi_i)$, then from Proposition \ref{prop MN2 3.10} (iii),
\begin{align}
\mathbf{M}\big(\Phi_i(x)\big)&\leq \mathbf{M}\Big(\Phi_i(x)- \Phi_i\big(\alpha_i(x)\big)\Big)+ \mathbf{M}\Big(\Phi_i\big(\alpha_i(x)\big)\Big) \nonumber \\
&\leq C_0(M, m)\mathbf{f}(\phi_i)+ \mathbf{M}\Big(\phi_i\big(\alpha_i(x)\big)\Big) \label{equi 1.5.5}
\end{align}
where $\alpha_i(x)\in X(\tau_i)_0$, and $\alpha_i(x), x$ lie in a common cell of $X(\tau_i)$.

From the definition of $\mathscr{D}_{p, m}$, when $i$ is large enough, $\Phi_i\in \mathscr{P}_{p, m}$. And from Proposition \ref{prop MN2 3.10}, $\Phi_i: X\rightarrow \mathcal{Z}_n^0(M; \mathbf{M}; \mathbb{Z}_2)$ is a continuous map. We have 
\begin{align}
\Phi_i\in \mathscr{P}_{p, m}^{\mathbf{M}} \ , \quad \quad \quad \forall i>> 1 \nonumber 
\end{align} 

Then by (\ref{equi 1.5.5}), the definition of $\omega_{p, m}^{\mathbf{M}}(M)$ and $\lim\limits_{i\rightarrow \infty} \mathbf{f}(\phi_i)\leq \lim\limits_{i\rightarrow \infty}\delta_i= 0$, 
\begin{align}
\omega_{p, m}^{\mathbf{M}}(M)&\leq \varlimsup_{i\rightarrow \infty} \max_{x\in dmn(\Phi_i)} \mathbf{M}\big(\Phi_i(x)\big) \leq \varlimsup_{i\rightarrow \infty} \max_{y\in dmn(\phi_i)} \mathbf{M}\big(\phi_i(y)\big)= \mathbf{L}(S_{\epsilon}) \label{equi 1.8}
\end{align}

From (\ref{equi 1.8}) and (\ref{equi 1.2}), we have
\begin{align}
\omega_{p, m}^{\mathbf{M}}(M)\leq \omega_{p, m}(M) \label{equi 1.9}
\end{align}

From (\ref{equi 1.1}) and (\ref{equi 1.9}), the first conclusion is proved. By (\ref{def p-width 1}) and (\ref{def p-M-width}), we get our second conclusion. 
}
\qed

\begin{lemma}\label{lem when P_(p, m) is non-empty}
{If $m\geq 2p+ 1$, then $\mathscr{P}_{p, m}^{\mathbf{M}}\neq \emptyset$.
}
\end{lemma}

\pf
{By \cite[Theorem $5.1$]{MN2}, there exists $\Phi: \mathbb{RP}^p\rightarrow \mathcal{Z}_n^0(M; \mathbb{Z}_2)$ such that $\Phi\in \mathscr{P}_{p}$. By Whitney's embedding theorem, $\mathbb{RP}^{p}$ can be embedded into $\mathbb{R}^{2p+ 1}$, hence we have $\Phi\in \mathscr{P}_{p, 2p+ 1}$. From Proposition \ref{prop width equation},
\begin{align}
\omega_{p, 2p+ 1}^{\mathbf{M}}(M)= \omega_{p, 2p+ 1}(M)< \infty \nonumber 
\end{align}
then $\mathscr{P}_{p, 2p+ 1}^{\mathbf{M}}\neq \emptyset$. If $m\geq 2p+ 1$, the conclusion is obtained by $\mathscr{P}_{p, 2p+ 1}^{\mathbf{M}}\subset \mathscr{P}_{p, m}^{\mathbf{M}}$.
}
\qed

\section{The existence of stationary and $\mathbb{Z}_2$ almost minimizing varifolds}\label{SECTION 4}

\begin{definition}\label{def discrete p-seqeunce}
{We define $\check{\mathscr{D}}_{p, m}$ as the set of $S= \{\phi_i\}_{i\in N}$,
\begin{align}
\phi_i: X^{(i)}(k_i)_0\rightarrow \mathcal{Z}_n(M; \mathbf{M}; \mathbb{Z}_2), \quad \quad \quad \quad \lim\limits_{i\rightarrow \infty}\mathbf{f}(\phi_i)= 0 \nonumber  
\end{align}
where $X^{(i)}$ is a cubical subcomplex of $I^m$, and for $\delta_0= \delta_0(M)$ as in Proposition \ref{prop MN2 3.10},
\begin{align}
\mathbf{f}(\phi_i)< \delta_0 \quad and \quad  \Phi_i= \mathscr{A}(\phi_i)\in \mathscr{P}_{p, m}^{\mathbf{M}} \ , \quad \quad \quad \quad \forall\ i\in \mathbb{N} \nonumber 
\end{align}
such $S$ is called a \textbf{$(p, m)$-homology sequence}. The \textbf{width} of $S$ is similarly defined as $\mathbf{L}(S)= \varlimsup_{i\rightarrow \infty}\max_{x\in dmn(\phi_i)}\mathbf{M}\left(\phi_i(x)\right)$.
}
\end{definition}

Similar as the proof of Lemma \ref{lem MN2's lem 4.7}, we have $\omega^{\mathbf{M}}_{p, m}(M)= \inf_{S\in \check{\mathscr{D}}_{p, m}} \mathbf{L}(S)$.

\begin{definition}\label{def p-critical sequence and set}
{$S= \{\phi_i\}_{i\in \mathbb{N}}$ is a \textbf{$(p, m)$-critical sequence} if
\begin{align}
S\in \check{\mathscr{D}}_{p, m} \quad \quad and \quad \quad \mathbf{L}(S)= \omega^{\mathbf{M}}_{p, m}(M) \nonumber 
\end{align}
For each $S= \{\phi_i\}_{i\in \mathbb{N}}\in \check{\mathscr{D}}_{p, m}$, define the \textbf{$(p, m)$-critical set} of $S$ as 
\begin{align}
\mathcal{C}_{p, m}(S)= \mathcal{K}_{p, m}(S)\cap \left\{V\in \mathcal{V}_n(M): \|V\|(M)= \mathbf{L}(S)\right\} \nonumber 
\end{align}
where 
\begin{align}
\mathcal{K}_{p, m}(S)&= \big\{V\in \mathcal{V}_n(M): \ V= \lim_{j\rightarrow \infty}\left|\phi_{i_j}(x_j)\right| \ for \ some \ increasing\ sequence \nonumber \\
&\quad  \{i_j\}_{j= \mathbb{N}}\ and \ some \ x_j\in dmn(\phi_{i_j})\big\} \nonumber  
\end{align}
}
\end{definition}

\begin{remark}\label{rem p-critical seq exists}
{For any $S\in \check{\mathscr{D}}_{p, m}$, it is easy to see that the critical set $\mathcal{C}_{p, m}(S)\neq \emptyset$ and is compact. If $\mathscr{P}_{p, m}^{\mathbf{M}}\neq \emptyset$, then $\check{\mathscr{D}}_{p, m}\neq \emptyset$ and there exists at least one $(p, m)$-critical sequence.
}
\end{remark}

Let $\mathscr{X}(M)$ denote the set of smooth vector fields on $M$, given $\chi\in \mathscr{X}(M)$, and let $\psi$ be the differential isotopy generated by $\chi$, i.e. $\frac{\partial}{\partial t}\psi= \chi(\psi)$. The \textbf{first variation} of $V$ with respect to $\chi$ in $M$ is defined as $[\delta V](\chi)= \frac{d}{dt} \big(\|\psi(t, \cdot)_{\#}V\|(M)\big)\big|_{t= 0}$. We say that the varifold $V$ is \textbf{stationary} in $M$ if 
\begin{align}
[\delta V](\chi)= 0\ , \quad \quad \quad \quad \quad \forall \chi\in \mathscr{X}(M) \nonumber 
\end{align}

\begin{prop}[Existence of stationary critical sets]\label{prop existence of stationary varifolds}
{For any $m, p\in \mathbb{N}$ with $m\geq 2p+ 1$, there exists a $(p, m)$-critical sequence $S\in \check{\mathscr{D}}_{p, m}$ such that $\mathcal{C}_{p, m}(S)\neq \emptyset$ and each $V\in \mathcal{C}_{p, m}(S)$ is stationary in $M$.
}
\end{prop}

\pf
{By Lemma \ref{lem when P_(p, m) is non-empty} and Remark \ref{rem p-critical seq exists}, there exists at least one $(p, m)$-critical sequence when $m\geq 2p+ 1$. Let $S^{*}= \{\varphi_i^{*}\}$ be a $(p, m)$-critical sequence and $C= \max_{i\in \mathbb{N}}\max_{x\in dmn(\varphi_i^{*})} \mathbf{M}\left(\varphi_i^{*}(x)\right)< \infty$. 

The rest proof of the conclusion is similar as \cite[Theorem $4.3$]{Pitts} (also see \cite[Section $15$]{MN1}), we obtained the existence of stationary varifolds by the classical pulling tight method. The key observation here is that the pulling tight method works well while replacing homotopy sequence by homology sequence, even without restricting the same value of $m$. The only difference with the former proofs is that the domains of $\varphi_i^*$ may be different. 
}
\qed

The following two combinatorics lemmas are needed to prove Proposition \ref{prop existence of almost minimizing varifolds}.

\begin{lemma}\label{lem Pitts 4.8}
{Suppose $\mathbf{A}= \big\{\bar{A}_{p_i}(s_{ij}, r_{ij}): i= 1, \cdots, k, j= 1, \cdots, km\big\}$ is a collection of closed annuli in $M$ such that 
\begin{align}
r_{ij}> s_{ij}> 3r_{i(j+ 1)}> 3s_{i(j+ 1)} \ , \quad \quad \quad \quad i= 1, \cdots, k\ , \ j= 1, \cdots, km-1 \label{ineq assumption}
\end{align}
Then there exists a disjointed collection $\mathbf{C}\subset \mathbf{A}$ such that 
\begin{align}
m= card \Big[\mathbf{C}\cap \big\{\bar{A}_{p_i}(s_{ij}, r_{ij}): j= 1, \cdots, km\big\}\Big] \ , \quad \quad \quad \quad i= 1, \cdots, k \nonumber 
\end{align}
}
\end{lemma}

\pf
{\cite[Lemma $4.8$]{Pitts}.
}
\qed

If $\beta= \bar{A}_p(s_1, s_2)$ is a closed annulus in $M^{n+ 1}$, for any $r> 0$, we define the map $\mu(r)$ as the following
\begin{align}
\mu(r)[\beta]= \big\{x\in M^{n+ 1}: rs_1\leq d_M(p, x)\leq rs_2\big\} \nonumber 
\end{align}

\begin{lemma}\label{lem Pitts 4.9}
{For each $\sigma\in I(m, k)$, let $\mathbf{A}(\sigma)$ be a disjointed collection of $(3^m)^{3^m}$ concentric closed annuli in $M^{n+ 1}$ such that
\begin{align}
a_1\cap \mu(r)[a_2]= \emptyset\ , \quad \quad \quad \quad \forall a_1, a_2\in \mathbf{A}(\sigma), a_1\neq a_2, 1\leq r\leq 3 \nonumber 
\end{align} 
Then there exists a map $\alpha$ defined on $I(m, k)$ such that 
\begin{align}
&\alpha(\sigma)\in \mathbf{A}(\sigma) \ , \quad \quad \quad \quad \forall \sigma\in I(m, k) \nonumber \\
&\alpha(\sigma)\cap \alpha(\tau)= \emptyset\ ,  \quad \quad if\ \sigma, \tau\in I(m, k)\ ,\ \sigma\neq \tau \quad and \quad \sigma, \tau\ are\ faces\ of \ Y\in I(m, k) \nonumber 
\end{align}
}
\end{lemma}

\pf
{\cite[Proposition $4.9$]{Pitts}.
}
\qed

For $T\in \mathcal{I}_k(M^{n+ 1})$ and $f: M^{n+ 1}\rightarrow \mathbb{R}$ is Lipschitz, we define \textbf{the slice of $T$ by $f$ at $r$} as $\langle T, f, r\rangle= \partial \Big(T\mres \big\{x: f(x)\leq r\big\}\Big)- (\partial T)\mres \big\{x: f(x)\leq r\big\}$. We have the following technical lemma:
\begin{lemma}\label{lem slice ineq}
{For $T\in \mathcal{I}_{n+ 1}(M^{n+ 1})$ and $f: M^{n+ 1}\rightarrow \mathbb{R}$ is Lipschitz, then 
\begin{align}
&\langle T, f, t\rangle\in \mathcal{I}_{n}(M^{n+ 1}) \ , \quad \quad \quad \quad for\ \mathscr{L}^1\ almost \ all \ t\in \mathbb{R} \nonumber \\
&\int_a^b \mathbf{M}\langle T, f, t\rangle dt \leq Lip(f)\cdot \mu_T\Big(\big\{x: a< f(x)< b\big\}\Big)\ , \quad \quad \quad \quad -\infty\leq a< b\leq \infty \nonumber 
\end{align}
}
\end{lemma}

\pf
{See \cite[Lemma $28.1$ and Lemma $28.5$]{Simon-book}.
}
\qed

For $x= (x_1, \cdots, y_m), y= (y_1, \cdots, y_m)\in I(m, j)_0$, we define 
\begin{align}
\mathfrak{d}(x, y)= 3^j\sum_{i= 1}^m |x_i- y_i| \nonumber 
\end{align}

\begin{definition}\label{def almost minimizing}
{For open subset $U\subset \subset M$, we say that $V\in \mathcal{V}_n(M)$ is \textbf{$\mathbb{Z}_2$ almost minimizing} in $U$ if and only if for each $\epsilon> 0$, there exists $\delta> 0$ and $T\in \mathcal{Z}_n(M, M- U; \mathbb{Z}_2)$ such that 
\begin{enumerate}
\item[(1)] $\mathbf{F}(V, |T|)< \epsilon$.
\item[(2)] If $T= T_0, T_1, \cdots, T_m\in \mathcal{Z}_n(M, M- U; \mathbb{Z}_2)$ and for   $i= 1, \cdots, m$,
\begin{align}
spt(T- T_i)\subset U\ , \quad \quad \mathbf{M}(T_i, T_{i- 1})\leq \delta\ , \quad \quad \mathbf{M}(T_i)\leq \mathbf{M}(T)+ \delta \ ,  \nonumber 
\end{align}
then $\mathbf{M}(T_m)> \mathbf{M}(T)- \epsilon$.
\end{enumerate}
}
\end{definition}

\begin{remark}\label{rem equ of almost minimizing def}
{By \cite[Theorem $3.9$]{Pitts}, the above definition is equivalent to the definition given in \cite[$3.1$]{Pitts}.
}
\end{remark}

For $p\in M$, $s> r>0$, define 
\begin{align}
B_p(r)&= \{x\in M: d_{M}(x, p)< r\} \nonumber \\
A_{p}(r, s)&= \{x\in M: r< d_{M}(x, p)< s\} \nonumber 
\end{align}

\begin{definition}\label{def almost minizing in annuli}
{A varifold $V\in \mathcal{V}_n(M)$ is \textbf{$\mathbb{Z}_2$ almost minimizing in annuli} if for each $z\in M$, there exists $r= r(z)> 0$ such that $V$ is $\mathbb{Z}_2$ almost minimizing in $A_z(s, r)$ for all $0< s< r$. 
}
\end{definition}

\begin{prop}\label{prop existence of almost minimizing varifolds}
{For $m, p\in \mathbb{N}$ with $m\geq 2p+ 1$, there exists $V\in \mathcal{V}_n(M)$ such that 
\begin{enumerate}
\item[$(1)$] $V$ is stationary in $M$ and $\|V\|(M)= \omega_{p, m}^{\mathbf{M}}(M)$.
\item[$(2)$] $V$ is $\mathbb{Z}_2$ almost minimizing in annuli.
\end{enumerate}
}
\end{prop}

\begin{remark}
{We follow the main line of argument in \cite[$4.10$]{Pitts}, the crucial difference is that the domain $X^{(i)}_0$ of every $\varphi_i$ in the $(p, m)$-critical sequence $S= \{\varphi_i\}$ possibly belongs to different dimensional cubical subcomplexes $X^{(i)}$ of $I^m$, the key observation is that all argument except the part involving the combinatorics lemmas of \cite{Pitts}, proceeds without essential problem. To use the combinatorics lemmas of \cite{Pitts}, the uniform upper bound of $m$ is enough.

}
\end{remark}

\pf
{\textbf{Step $1$. Choose $i$ large enough.}

From Proposition \ref{prop existence of stationary varifolds}, we can obtain a $(p, m)$-critical sequence $S= \{\varphi_i\}$, such that if $V\in \mathcal{C}_{p, m}(S)$, then $\|V\|(M)= \omega_{p, m}^{\mathbf{M}}(M)$ and $V$ is stationary in $M$. 

We assume that no $V\in \mathcal{C}_{p, m}(S)$ satisfies conclusion $(2)$ of Proposition \ref{prop existence of almost minimizing varifolds} and obtain the contradiction by constructing a $(p, m)$-homology sequence $S^{*}$ with $\mathbf{L}(S^{*})< \mathbf{L}(S)$. 

We abbreviate $L= 2^m$, $c= (3^m)^{3^m}$. By our assumption, for each $V\in \mathcal{C}_{p, m}(S)$, there exists $z\in spt (\|V\|)$ and positive numbers $r_1, \cdots, r_c, s_1, \cdots, s_c$ such that 
\begin{align}
r_c- 2s_c> 0 \quad and \quad r_k- 2s_k> 3(r_{k+ 1}+ 2s_{k+ 1})\ , \quad \quad \quad \quad k= 1, \cdots, c- 1 \nonumber 
\end{align}
$V$ is not $\mathbb{Z}_2$ almost minimizing in $A_z(r_k- s_k, r_k+ s_k)$, $k= 1, \cdots, c$. We set 
\begin{align}
A_k(V)&= A_z(r_k- 2s_k, r_k+ 2s_k)\ , \quad \quad \quad \quad a_k(V)= A_z(r_k- s_k, r_k+ s_k) \nonumber \\
s(V)&= \min\{s_1, \cdots, s_c\} \nonumber 
\end{align}

For each $V\in \mathcal{C}_{p, m}(S)$, there exists $\epsilon(V, c)> 0$ such that for each $j= 1, \cdots, c$, each $\delta> 0$, and each $T\in \mathcal{Z}_n(M; \mathbf{M}; \mathbb{Z}_2)$ with $\mathbf{F}\left(V, |T|\right)< \epsilon(V, c)$, there exists a sequence $T= T_1, \cdots, T_q\in \mathcal{Z}_n(M; \mathbf{M}; \mathbb{Z}_2)$, where $q$ depends on $T$ and $\delta$, such that 
\begin{align}
&\bigcup_k spt(T_k- T) \subset a_j(V)\ , \quad \quad \quad \quad \sup_k \mathbf{M}(T_k- T_{k- 1})\leq \delta \label{part 2} \\
&\sup_k \mathbf{M}(T_k)\leq \mathbf{M}(T)+ \delta \ , \quad \quad \quad \quad \mathbf{M}(T_q)< \mathbf{M}(T)- \epsilon (V, c) \label{part 2.1}
\end{align}
Because $c= (3^m)^{3^m}$ is fixed in the proof, we will use $\epsilon(V)$ instead of $\epsilon(V, c)$ in the following of the proof.

By $\mathcal{C}_{p, m}(S)$ is compact, there exists $V_1, \cdots, V_{\ell}$ in $\mathcal{C}_{p, m}(S)$ such that 
\begin{align}
\mathcal{C}_{p, m}(S)\subset \bigcup_{k= 1}^{\ell} \left\{V\in \mathcal{V}_n(M): \mathbf{F}(V, V_k)< \frac{1}{4} \epsilon(V_k) \right\} \nonumber 
\end{align}

We define positive numbers $\epsilon_1, \epsilon_2, s$ as follows:
\begin{align}
\epsilon_1= \min_{k= 1, \cdots, \ell} \epsilon (V_k) \ , \quad \quad \quad \quad \epsilon_2= \min\left\{\frac{1}{2}\epsilon_1, \frac{1}{4}\tilde{\epsilon}\right\}\ , \quad \quad \quad \quad s= \inf_{j= 1, \cdots, \ell} s(V_j)  \nonumber 
\end{align}
where 
\begin{align}
\tilde{\epsilon}&= \sup\big\{\epsilon: \mathbf{F}(V, V_k)< \frac{1}{4}\epsilon(V_k)\ for \ some \ k= 1, \cdots, \ell, \nonumber \\
&\quad \quad \quad whenever\ V\in \mathcal{K}_{p, m}(S)\ and \ \|V\|(M)\geq \mathbf{L}(S)- 2\epsilon \big\} \nonumber 
\end{align}

We assume that $dmn(\varphi_i)= X^{(i)}(n_i)_0$, $i= 1, 2, \cdots$. We choose a positive integer $N$ such that if $i\geq N$, the following properties hold:
\begin{enumerate}
\item[$(N.1)$] Either $\mathbf{M}(\varphi_i(x))< \mathbf{L}(S)- 2\epsilon_2$ or 
\begin{align}
\mathbf{M}(\varphi_i(x))\geq \mathbf{L}(S)- 2\epsilon_2 \quad \quad and \quad \quad 
\mathbf{F}\left(|\varphi_i(x)|, V_j\right)< \frac{1}{2}\epsilon (V_j) \ for \ some \ j \nonumber 
\end{align}
\item[$(N.2)$] $\mathbf{f}(\varphi_i)\leq \epsilon_2$, and by \cite[Corollary $1.14$]{Alm}, if $\sigma\in X^{(i)}(n_i)$, $\{x, y\}\subset \sigma_0$, then there exists $Q\in \mathcal{I}_{n+ 1}(M; \mathbb{Z}_2)$ with 
\begin{align}
\partial Q= \varphi_i(x)- \varphi_i(y) \quad \quad and \quad \quad \mathbf{M}(Q)\leq \mathbf{M}\left(\varphi_i(x)- \varphi_i(y)\right) \label{part 5(c)} 
\end{align}
\item[$(N.3)$] $\mathbf{f}(\varphi_i)< \frac{1}{2m}\epsilon_1$. 
\end{enumerate}

In the rest of the proof, we always assume $i\geq N$ and fix $i$.

\textbf{Step $2$. Construct discrete deformation from non $\mathbb{Z}_2$ almost minimizing.}

Whenever $x\in X^{(i)}(n_i)_0$ and $\mathbf{M}\left(\varphi_i(x)\right)\geq \mathbf{L}(S)- 2\epsilon_2$, choose $f_1(x)\in \{1, 2, \cdots, \ell\}$ such that $\mathbf{F}\left(|\varphi_i(x)|, V_{f_1(x)}\right)< \frac{1}{2}\epsilon\left(V_{f_1(x)}\right)$. Let $\delta_i= 2L^2 m\mathbf{f}(\varphi_i)\left(1+ 4(L- 1)s^{-1}\right)$. We have the following claim:
\begin{claim}\label{claim existence of sequence in non almost minimizing}
{There is a positive integer $N_1= N_1(i)$ such that whenever $\sigma\in X^{(i)}(n_i)$, $\mu= 1, \cdots, c$, $\{x_1, \cdots, x_{\ell}\}\subset \sigma_0$, $\min\limits_{j= 1, \cdots, \ell} \mathbf{M}\left(\varphi_i(x_j)\right)\geq \mathbf{L}(S)- 2\epsilon_2$, there exists a sequence $\{T_{j,q}\}_{j= 1, \cdots, \ell \atop q=1, \cdots, 3^{N_1}} \subset \mathcal{Z}_n(M; \mathbf{M}; \mathbb{Z}_2)$ such that for any $j, j'=1, \cdots, \ell$,
\begin{align}
&\varphi_i(x_j)= T_{j,1}\ , \quad \quad \quad \quad \bigcup_{q= 1}^{3^{N_1}} spt\left(T_{j, q}- T_{j, 1}\right)\subset A_{\mu}\left(V_{f_1(x_1)}\right)\ , \nonumber \\
&\mathbf{M}\left(T_{j, q}- T_{j, q- 1}\right)\leq L^{- 2}\delta_i\ , \quad \quad \quad \quad q= 2, 3, \cdots, 3^{N_1} \nonumber \\
&\mathbf{M}\left(T_{j, q}- T_{j', q}\right)\leq L^{- 2}\delta_i\ , \quad \quad \quad \quad q= 1, 2, \cdots, 3^{N_1} \nonumber \\
&\mathbf{M}\left(T_{j, q}\right)\leq \mathbf{M}\left(T_{j, 1}\right)+ L^{-2}\delta_i\ , \quad \quad \quad \quad q= 1, 2, \cdots, 3^{N_1} \nonumber \\
& \mathbf{M}(T_{j, 3^{N_1}})< \mathbf{M}(T_{j, 1})- \big(\epsilon_1- 2L^{-2}\delta_i\big) \nonumber 
\end{align}
}
\end{claim}

\pf
{Note $X^{(i)}$ is a cubical subcomplex of $I^m$, we in fact have $\ell\leq 2^m= L$. And for $j= 1, 2, \cdots, \ell$, from the definition of $\mathbf{f}$, we have 
\begin{align}
\mathbf{F}\left(|\varphi_i(x_j)|, V_{f_1(x_1)}\right)&\leq \mathbf{M}\left(\varphi_i(x_j)- \varphi_i(x_1)\right)+\mathbf{F}\left(|\varphi_i(x_1)|, V_{f_1(x_1)}\right) \nonumber \\
&\leq m\mathbf{f}(\varphi_i)+ \frac{1}{2}\epsilon\left(V_{f_1(x_1)}\right)< \frac{1}{2}\epsilon_1+ \frac{1}{2}\epsilon(V_{f_1(x_1)})\leq \epsilon\left(V_{f_1(x_1)}\right)  \nonumber 
\end{align}

As (\ref{part 2}) and (\ref{part 2.1}), we can get a sequence $T_1\vcentcolon =\varphi_i(x_1), T_2, \cdots, T_q, T_{q+ 1}\in \mathcal{Z}_n(M; \mathbf{M}; \mathbb{Z}_2)$, where $q+ 1= 3^{N_1}$, such that $T_q= T_{q+ 1}$ and 
\begin{align}
\bigcup_k spt(T_k- T_1) \subset a_{\mu}(V)\ , \quad \quad \quad \quad \sup_k \mathbf{M}(T_k- T_{k- 1})\leq \mathbf{f}(\varphi_i) \nonumber \\
\sup_k \mathbf{M}(T_k)\leq \mathbf{M}(T_1)+ \mathbf{f}(\varphi_i) \ , \quad \quad \quad \quad \mathbf{M}(T_q)< \mathbf{M}(T_1)- \epsilon (V) \nonumber 
\end{align}

If $a_{\mu}\left(V_{f_1(x_1)}\right)= A_p(r_{\mu}- s_{\mu}, r_{\mu}+ s_{\mu})$, then by (\ref{part 5(c)}), there exist $Q_2, \cdots, Q_\ell\in \mathcal{I}_{n+ 1}(M; \mathbb{Z}_2)$ such that for each $j= 2, \cdots, \ell$, 
\begin{align}
\partial Q_j&= \varphi_i(x_j)- \varphi_i(x_1) \nonumber \\
\mathbf{M}(Q_j)&\leq \mathbf{M}\left(\varphi_i(x_j)- \varphi_i(x_1)\right)\leq m\mathbf{f}(\varphi_i) \nonumber 
\end{align}

Writing $u(x)= d_M(x, p)$, $x\in M$, and from Lemma \ref{lem slice ineq},
\begin{align}
\int_a^b \mathbf{M}\langle Q_j, u, t\rangle dt \leq \mathbf{M}(Q_j)\ , \quad \quad \quad \quad a< b \nonumber 
\end{align}
we infer that there exist $t_1, t_2\in \mathbb{R}$ such that 
\begin{align}
&r_{\mu}+ s_{\mu}< t_1< r_{\mu}+ 2s_{\mu}\ , \quad \quad r_{\mu}- 2s_{\mu}< t_2< r_{\mu}- s_{\mu}\ , \nonumber \\
&\langle Q_j, u, t_k\rangle\in \mathcal{I}_n(M; \mathbb{Z}_2)\ , \quad \quad \quad \quad k= 1, 2, \quad  j= 2, \cdots, \ell \nonumber \\
&\mathbf{M}\langle Q_j, u, t_k\rangle\leq (L- 1)m\mathbf{f}(\varphi_i)s_{\mu}^{-1}\ ,  \quad \quad \quad \quad k= 1, 2, \quad  j= 2, \cdots, \ell \nonumber
\end{align}

Define $T_{1, l}= T_l$ for any $l= 1, \cdots, q+1$. For each $j= 2, \cdots, \ell$, we define $T_{j, 1}= \varphi_i(x_j)$ and if $2\leq k\leq q+ 1$,
\begin{align}
T_{j, k}= \varphi_i(x_j)\mres \left(\left[M^{n+ 1}- \bar{B}_p(t_1)\right] \cup B_p(t_2)\right)+ \langle Q_j, u, t_2\rangle- \langle Q_j, u, t_1\rangle + T_{k- 1} \mres A_p(t_2, t_1)  \nonumber 
\end{align}

Then we can verify the conclusion using the above $T_{j, k}$.
}
\qed

\textbf{Step $3$. Choose the suitable annuli on $M^{n+ 1}$.}

We choose a positive integer $N_2= N_2(i)\geq n_i$ and define 
\begin{align}
f_2: X^{(i)}(N_2)_0\rightarrow \{0, 1, \cdots, 3^{N_1}\} \nonumber 
\end{align}
such that 
\begin{enumerate}
\item[$(f_2.1)$] 
\begin{equation}\nonumber 
f_2(x)= \left\{
\begin{array}{rl}
& 0 \ , \quad \quad \quad \quad \quad \quad \quad \quad \quad \quad \quad\quad \quad \mathbf{M}\left(\varphi_i\circ \mathbf{n}(N_2, n_i)(x)\right)< \mathbf{L}(S)- 2\epsilon_2 \\
& 3^{N_1} \ , \quad \quad \quad \quad \quad \quad \quad \quad \quad \quad \quad\quad \quad \mathbf{M}\left(\varphi_i\circ \mathbf{n}(N_2, n_i)(x)\right)\geq  \mathbf{L}(S)- \epsilon_2
\end{array} \right.
\end{equation}
\item[$(f_2.2)$] If $\sigma\in X^{(i)}(N_2)$, $\{x, y\}\subset \sigma_0$, then $\left|f_2(x)- f_2(y)\right|\leq 1$, and if furthermore $f_2(x)> 0$, 
\begin{align}
\mathbf{M}\big(\varphi_i\circ \mathbf{n}(N_2, n_i)(y)\big)\geq \mathbf{L}(S)- 2\epsilon_2  \label{4.28.1.1}
\end{align}
\end{enumerate}

We define 
\begin{align}
f_4&: \left\{\tau\in X^{(i)}(N_2): \sup_{x\in \tau_0}f_2(x)> 0 \right\}\rightarrow \{1, \cdots, \ell \} \nonumber \\
f_4&= f_1\circ \mathbf{n}(N_2, n_i)\circ f_3 \nonumber 
\end{align}
where $f_3: X^{(i)}(N_2)\rightarrow X^{(i)}(N_2)_0$ is any map satisfying $f_3(\sigma)\in \sigma_0$.

Let $2^M$ denote the family of all subsets of $M$, and one can define 
\begin{align}
f_5: \left\{\tau\in X^{(i)}(N_2): \sup_{x\in \tau_0}f_2(x)> 0 \right\}\rightarrow 2^{M}\nonumber 
\end{align}
such that 
\begin{enumerate}
\item[$(f_5.1)$] $f_5(\tau)= A_j\left(V_{f_4(\tau)}\right)$ for some $j\in \{1, \cdots, c\}$.
\item[$(f_5.2)$] If $\sigma\neq \tau$ and $\sigma, \tau$ are faces (not necessarily of the same dimension) of a common cell in $X^{(i)}(N_2)$, then $dist\left(f_5(\sigma), f_5(\tau)\right)> 0$. 
\end{enumerate}
These choices are possible by Lemma \ref{lem Pitts 4.9}. 

\textbf{Step $4$. Match the deformation with the annuli.}

Let $\tau\in X^{(i)}(N_2)$ and $\mu\in \{0, 1, \cdots, 3^{N_1}\}$. For each $x\in \tau_0$, we define $f_6(\tau, x, \mu)\in \mathcal{Z}_n(M; \mathbf{M}; \mathbb{Z}_2)$ as follows.
\begin{enumerate}
\item[$(f_6.1)$] If $\sup\limits_{x\in \tau_0} f_2(x)= 0$, then $f_6(\tau, x, \mu)= \varphi_i\circ \mathbf{n}(N_2, n_i)(x)$.
\item[$(f_6.2)$] If $\sup\limits_{x\in \tau_0} f_2(x)> 0$, then by (\ref{4.28.1.1}) and Claim \ref{claim existence of sequence in non almost minimizing}, for each $y\in \tau_0$, there exists a sequence 
\begin{align}
T_y(1)= \varphi_i\circ \mathbf{n}(N_2, n_i)(y), T_y(2), \cdots, T_y(3^{N_1})\in \mathcal{Z}_n(M; \mathbf{M}; \mathbb{Z}_2) \nonumber 
\end{align}
such that 
\begin{align}
& spt\left(T_y(1)- T_y(j)\right)\subset f_5(\tau)\ , \quad \quad \quad \quad j= 2, 3, \cdots, 3^{N_1} \nonumber \\
&\mathbf{M}\left(T_y(j)- T_y(j- 1)\right)\leq L^{- 2}\delta_i\ , \quad \quad \quad \quad j= 2, 3, \cdots, 3^{N_1} \nonumber \\
&\mathbf{M}\left(T_y(j)- T_z(j)\right)\leq L^{- 2}\delta_i\ , \quad \quad \quad \quad j= 1, 2, \cdots, 3^{N_1}\ , \ z\in \tau_0 \nonumber \\
&\mathbf{M}\left(T_y(j)\right)\leq \mathbf{M}\left(T_y(1)\right)+ L^{-2}\delta_i\ , \quad \quad \quad \quad j= 1, 2, \cdots, 3^{N_1} \nonumber \\
& \mathbf{M}\left(T_y(3^{N_1})\right)< \mathbf{M}\left(T_y(1)\right)- (\epsilon_1- 2L^{-2} \delta_i) \nonumber  
\end{align}

We set 
\begin{equation}\nonumber 
f_6(\tau, x, \mu)= \left\{
\begin{array}{rl}
& \varphi_i\circ \mathbf{n}(N_2, n_i)(x) \ , \quad \quad \quad \quad \quad \quad \quad \quad \quad \quad \quad\quad \quad \mu= 0 \\
& T_x(\mu) \ , \quad \quad \quad \quad \quad \quad \quad \quad \quad \quad \quad\quad \quad 1\leq \mu\leq \min\{3^{N_1}, f_2(x)\} \\
& T_x\left(\min\{3^{N_1}, f_2(x)\}\right)\ , \quad \quad \quad \quad \quad \quad \quad  \quad \quad \min\{3^{N_1}, f_2(x)\}\leq \mu\leq 3^{N_1} 
\end{array} \right.
\end{equation}
\end{enumerate}

\textbf{Step $5$. The choice of parameters in the homotopy.}

Let $N_3= N_3(i)= N_1+ N_2+ 2$. We define $f_7: X^{(i)}(N_3)_0\rightarrow X^{(i)}(N_2)$ so that for each $x\in X^{(i)}(N_3)_0$, $f_7(x)$ is the unique cell of least dimension in $X^{(i)}(N_2)$ containing $x$. We define $f_8(x, \tau)= \max\{0, 3^{N_1}- \gamma\}$, which is a map from $\big[X^{(i)}(N_3)_0\times X^{(i)}(N_2)\big] \cap \big\{(x, \tau): \tau \ is \ a \ face \ of \ f_7(x)\ and\ \mathbf{n}(N_3, N_2)(x)\in \tau_0\big\}$ to $\{0, 1, \cdots, 3^{N_1}\}$, where $\gamma$ is a number to be determined below.

\begin{enumerate}
\item[$(\gamma.1)$] If $\tau\in X^{(i)}(N_2)_0$, we set 
\begin{equation}\nonumber 
\gamma= \left\{
\begin{array}{rl}
& 0 \ , \quad \quad \quad \quad \quad \quad \quad \quad \quad \quad \quad\quad \quad |x- \tau|\leq 3^{-N_2- 1} \\
& \inf\left\{\mathfrak{d}(x, y): y\in X^{(i)}(N_3)_0\ and \ |y- \tau|\leq 3^{- 1- N_2}\right\} \ , \quad \quad  \quad |x- \tau|> 3^{-N_2- 1}
\end{array} \right.
\end{equation}
\item[$(\gamma.2)$] If $\tau$ is a $j$-cell , $j\geq 1$, let $x^{*}= (x_1^*, \cdots, x_m^*)$ be the unique element of $\tau_0$ with the property that 
\begin{align}
\mathfrak{d}\left([0], x^{*}\right)= \inf_{y\in \tau_0}\mathfrak{d}\left([0], y\right) \nonumber 
\end{align}
and let $\mathfrak{S}$ be the set of all $\omega= (\omega_1, \cdots, \omega_m)\in X^{(i)}(N_3)_0$ with the property:
\begin{equation}\nonumber 
\left\{
\begin{array}{rl}
& 3^{- 1- N_2}\leq |\omega_k- x_k^{*}|\leq 2\cdot 3^{- 1- N_2}\ , \quad \quad \quad \quad if\ |z_k- y_k|\neq 0\ for\ some \ \{z, y\}\subset \tau_0 \nonumber  \\
& |\omega_k- x_k^{*}|\leq 3^{- 1- N_2}\ ,\quad \quad \quad \quad for\ other \ k \nonumber 
\end{array} \right.
\end{equation}
we set $\gamma= \inf_{y\in \mathfrak{S}} \mathfrak{d}(x, y)$.
\end{enumerate}

\textbf{Step $6$. Construct the whole homotopy between $\{\varphi_i\}$ and $\{\varphi_i^{*}\}$.}

Now we define the discrete homotopy 
\begin{align}
\psi_i: I(1, N_3)_0\times X^{(i)}(N_3)_0\rightarrow \mathcal{Z}_n(M; \mathbf{M}; \mathbb{Z}_2) \nonumber 
\end{align}
as follows. 

\begin{enumerate}
\item[$(\psi_i.1)$] If $f_2\left(\mathbf{n}(N_3, N_2)(x)\right)= 0$, we set $\psi_i(j, x)= \varphi_i\circ \mathbf{n}(N_3, n_i)(x)$. 
\item[$(\psi_i.2)$] If $f_2\left(\mathbf{n}(N_3, N_2)(x)\right)> 0$, set $\xi(x, \tau)= \min\big\{f_2\circ \mathbf{n}(N_3, N_2)(x), f_8(x, \tau)\big\}$ and  
\begin{equation}\nonumber 
\left\{
\begin{array}{rl}
& \psi_i(j, x)\mres Z= \varphi_i\circ \mathbf{n}(N_3, n_i)(x)\mres Z \ , \quad \quad \quad \quad 0\leq j\cdot 3^{N_3}\leq 3^{N_1} \nonumber  \\
& \psi_i(j, x)\mres f_5(\tau)= \varphi_i\circ \mathbf{n}(N_3, n_i)(x)\mres f_5(\tau)\ , \quad \quad \quad \quad j= 0 , 3^{- N_3} \nonumber  \\
& \psi_i(j, x)\mres f_5(\tau)= f_6\left(\tau, \mathbf{n}(N_3, N_2)(x), j\cdot 3^{N_3}\right)\mres f_5(\tau)\ , \quad \quad \quad \quad if\ 1\leq j\cdot 3^{N_3}\leq \xi(x, \tau)\nonumber\\
& \psi_i(j, x)\mres f_5(\tau)= f_6\left(\tau, \mathbf{n}(N_3, N_2)(x), \xi(x, \tau)\right)\mres f_5(\tau)\ , \quad\quad \quad \quad  if\  \xi(x, \tau)\leq j\cdot 3^{N_3}\leq 3^{N_1} \nonumber \\
&\psi_i(j, x)= \psi_i(3^{N_1}\cdot 3^{- N_3}, x)\ , \quad \quad \quad \quad 3^{N_1}\leq j\cdot 3^{N_3}\leq 3^{N_3} \nonumber 
\end{array} \right.
\end{equation}
where $Z= M- \bigcup\left\{f_5(\tau): \tau \ is \ a \ face \ of \ f_7(x)\ and \ \mathbf{n}(N_3, N_2)(x)\in \tau_0 \right\}$. 
\end{enumerate}
One can verify that $\mathbf{f}(\psi_i)\leq \delta_i$.

Finally we define 
\begin{equation}\nonumber 
\varphi_i^{*}= \left\{
\begin{array}{rl}
& \varphi_i(x) \ , \quad \quad \quad \quad x\in dmn(\varphi_i)\ , \ i=1, 2, \cdots, N- 1 \nonumber  \\
& \psi_i(1, x)\ , \quad \quad \quad \quad x\in X^{(i)}(N_3)_0\ , \ i= N, N+1, \cdots \nonumber  
\end{array} \right.
\end{equation}

One can verify that $S^{*}= \{\varphi_i^{*}\}\in \check{\mathscr{D}}_{p, m}$ and $\mathbf{L}(S^{*})\leq \mathbf{L}(S)- \epsilon_2$, it is the contradiction to the fact that $S$ is the $(p, m)$-critical sequence.
}
\qed

{\it \textbf{Proof of Theorem \ref{thm main result}}}:~
{It follows from Proposition \ref{prop existence of almost minimizing varifolds}, that there exists $V\in \mathcal{IV}_n(M)$, which is $\mathbb{Z}_2$ almost minimizing in annuli and stationary in $M$, furthermore $\|V\|(M)= \omega_{p, m}^{\mathbf{M}}(M)$. Then apply Proposition \ref{prop width equation}, we get 
\begin{align}
\|V\|(M)= \omega_{p, m}(M) \nonumber 
\end{align}

Recall \cite[$3.13$]{Pitts} as the following:
\begin{lemma}\nonumber 
{Let $V\in \mathcal{V}_n(M)$, if $V$ is $\mathbb{Z}_2$ almost minimizing in annuli and $V$ is stationary in $M$, then $V\in \mathcal{IV}_n(M)$.
}
\end{lemma}
Apply it to $V$, we obtain that $V\in \mathcal{IV}_n(M^{n+ 1})$. Also we know \cite[Theorem $2.11$]{MN2} as the following:
\begin{theorem}\nonumber 
{Suppose $n\leq 6$, let $V\in \mathcal{IV}_n(M^{n+ 1})$ be a nontrivial integral varifold that is both stationary in $M^{n+ 1}$ and $\mathbb{Z}_2$ almost minimizing in annuli. Then $V$ is the varifold of a smooth, closed, embedded minimal hypersurface, with possible multiplicities.
}
\end{theorem}
The conclusion finally follows from the above theorem.
}
\qed

\appendix

\section*{Acknowledgments}
The author would like to thank Martin Li for suggestion, Zhiqin Lu, Fernando Marques and William Minicozzi for comments. He is indebted to Allan Hatcher and Leon Simon for help on algebraic topology and geometric measure theory respectively, Larry Guth and Jiaping Wang for encouragement. Especially, he is deeply grateful to Xin Zhou for the conversation and his suggestions. Finally, we thank the referee for the suggestion and comments, which are both helpful and insightful.

\end{document}